%%%%%%%%%%%
% typoref.tex. V : January 18, 2000.
% Author : Anthony PHAN
% Warning : syntaxe +- LaTeX
% Sources :
% T. Lachand--Robert, ``La Ma\^\i trise de \TeX'',
% R\'ef\'erences crois\'ees;
% latex.ltx's sources;
% and of course the \TeX book.
%%%%%%%%%%%%%%%%%%%%%%%%%%%%%%%%%%%%%%%%%%%%%%%%%%%%%
%

\catcode`\@=11

\magnification=1200
\baselineskip=14pt

\pretolerance=500    \tolerance=1000 \brokenpenalty=5000

\catcode`\;=\active
\def;{\relax\ifhmode\ifdim\lastskip>\z@
\unskip\fi\kern.2em\fi\string;}

\overfullrule=0mm

\catcode`\!=\active
\def!{\relax\ifhmode\ifdim\lastskip>\z@
\unskip\fi\kern.2em\fi\string!}

\catcode`\?=\active
\def?{\relax\ifhmode\ifdim\lastskip>\z@
\unskip\fi\kern.2em\fi\string?}

 \nonfrenchspacing

\newif\ifpagetitre            \pagetitretrue
\newtoks\hautpagetitre        \hautpagetitre={ }
\newtoks\baspagetitre         \baspagetitre={1}

\newtoks\auteurcourant        \auteurcourant={ }
\newtoks\titrecourant
\titrecourant={ }
\newtoks\hautpagegauche       \newtoks\hautpagedroite
\hautpagegauche={\hfill\sevenrm\the\auteurcourant\hfill}
\hautpagedroite={\hfill\sevenrm\the\titrecourant\hfill}

\newtoks\baspagegauche       \baspagegauche={\hfill\rm\folio\hfill}

\newtoks\baspagedroite       \baspagedroite={\hfill\rm\folio\hfill}

\headline={
\ifpagetitre\the\hautpagetitre
\global\pagetitrefalse
\else\ifodd\pageno\the\hautpagedroite
\else\the\hautpagegauche\fi\fi}

\footline={\ifpagetitre\the\baspagetitre
\global\pagetitrefalse
\else\ifodd\pageno\the\baspagedroite
\else\the\baspagegauche\fi\fi}

\def\date{\ {\the\day}\
\ifcase\month\or Janvier\or F\'evrier\or Mars\or Avril
\or Mai \or Juin\or Juillet\or Ao\^ut\or Septembre
\or Octobre\or Novembre\or D\'ecembre\fi\
{\the\year}}

\def\up#1{\raise 1ex\hbox{\sevenrm#1}}

\def\cqfd{\unskip\kern 6pt\penalty 500
\raise -2pt\hbox{\vrule\vbox to 10pt{\hrule width 4pt
\vfill\hrule}\vrule}\par\medskip}

\def\section#1{\vskip 7mm plus 20mm minus 1.5mm\penalty-50
\vskip 0mm plus -20mm minus 1.5mm\penalty-50
{\bf\noindent#1}\nobreak\smallskip}

\def\subsection#1{\medskip{\bf#1}\nobreak\smallskip}

\def\displaylinesno #1{\dspl@y\halign{
\hbox to\displaywidth{$\@lign\hfil\displaystyle##\hfil$}&
\llap{$##$}\crcr#1\crcr}}

\def\ldisplaylinesno #1{\dspl@y\halign{
\hbox to\displaywidth{$\@lign\hfil\displaystyle##\hfil$}&
\kern-\displaywidth\rlap{$##$}
\tabskip\displaywidth\crcr#1\crcr}}

\def\hfl#1#2{\smash{\mathop{\hbox to 12 mm{\rightarrowfill}}
\limits^{\scriptstyle#1}_{\scriptstyle#2}}}

%
% style (look at the behavior of \item dans \bibitem too,
% and at one ,\  in \re@dreferenceslist)
% Feel free to change: 	\bibn@me (title like ``R\'ef\'erences'')
%			\bibliographym@rk (general style)
%
\def\bibn@me{R\'ef\'erences}
\def\bibliographym@rk{\centerline{{\sc\bibn@me}}
	\sectionmark\section{\ignorespaces}{\unskip\bibn@me}
	\bigbreak\bgroup
	\ifx\ninepoint\undefined\relax\else\ninepoint\fi}
%
% Beware of the \bgroup: it will be closed by \endthebibliography
%
% \refsp@ce is the spacing command that appens between multiple
% references.
%
\let\refsp@ce=\
\let\bibleftm@rk=[
\let\bibrightm@rk=]
%
% if you want more space between brackets...
%\let\refsp@ce=\thinspace
%\def\bibleftm@rk{[\thinspace}
%\def\bibrightm@rk{\thinspace]}
%
% frenchy stuff
%
\def\numero{n\raise.82ex\hbox{$\fam0\scriptscriptstyle
o$}~\ignorespaces}
%
% new variables
%
\newcount\equationc@unt
\newcount\bibc@unt
\newif\ifref@changes\ref@changesfalse
\newif\ifpageref@changes\ref@changesfalse
\newif\ifbib@changes\bib@changesfalse
\newif\ifref@undefined\ref@undefinedfalse
\newif\ifpageref@undefined\ref@undefinedfalse
\newif\ifbib@undefined\bib@undefinedfalse
\newwrite\@auxout
%
% mark an equation
%
%\def\eqnum{\global\advance\equationc@unt by 1%
%\edef\lastref{\number\equationc@unt}%
%\eqno{(\lastref)}}
%
% One can reference anything, just copy the former macro
% and use it so: \machin \label{truc}
% In machin you would have defined \lastref by some number
% or any text.
%
% References macros
%
% The next macros are the core of \ref and \cite commands.
% Its first argument may be ref, pageref or bib.
%
% It is too tricky to be explained.
% (It is a bit recursive.)
% It allows using \cite or \ref or ...
% with arbitrary many arguments,
% for instance:
% \cite{knuth1,knuth2,ma pomme}
%
% First argument is always ref, pageref or bib.
%
\def\re@dreferences#1#2{{%
	\re@dreferenceslist{#1}#2,\undefined\@@}}
\def\re@dreferenceslist#1#2,#3\@@{\def\next{#2}%
	\expandafter\ifx\csname#1@@\meaning\next\endcsname\relax
	??\immediate\write16
	{Warning, #1-reference "\next" on page \the\pageno\space
	is undefined.}%
	\global\csname#1@undefinedtrue\endcsname
	\else\csname#1@@\meaning\next\endcsname\fi
	\ifx#3\undefined\relax
	\else,\refsp@ce\re@dreferenceslist{#1}#3\@@\fi}
%
% notice that the former ``,\refsp@ce'' will separate
% multiple arguments. But beware of spaces
% while defining a reference or calling for it!
%
% tricky thing: \newlabel has two arguments
% {labelname}{{\lastref}{\pageref}}
% The second argument is read as two arguments
% by \newl@bel. This was necessary to get
% a jobname.aux containing the same syntax
% LaTeX would produce and use.
%
\def\newlabel#1#2{{\def\next{#1}\newl@bel#2}}
\def\newl@bel#1#2{%
	\expandafter\xdef\csname ref@@\meaning\next\endcsname{#1}%
	\expandafter\xdef\csname pageref@@\meaning\next\endcsname{#2}}
\def\label#1{{%
	\toks0={#1}\message{ref(\lastref) \the\toks0,}%
	\ignorespaces\immediate\write\@auxout%
	{\noexpand\newlabel{\the\toks0}{{\lastref}{\the\pageno}}}%
	\def\next{#1}%
	\expandafter\ifx\csname ref@@\meaning\next\endcsname\lastref%
	\else\global\ref@changestrue\fi%
	\newlabel{#1}{{\lastref}{\the\pageno}}}}
\def\ref#1{\re@dreferences{ref}{#1}}
\def\pageref#1{\re@dreferences{pageref}{#1}}
%
% bibliography macros
%
\def\bibcite#1#2{{\def\next{#1}%
	\expandafter\xdef\csname bib@@\meaning\next\endcsname{#2}}}
\def\cite#1{\bibleftm@rk\re@dreferences{bib}{#1}\bibrightm@rk}
%
% The argument of \beginthebibliography
% is any sequence of numerals which will represent
% the maximum \item's length. If you have less than 9
% \bibitem's, this argument may be {any numeral}.
% if you have between 100 and 999 \bibitem's
% this argument may be {any three numerals},
% and so on.
%
\def\beginthebibliography#1{\bibliographym@rk
	\setbox0\hbox{\bibleftm@rk#1\bibrightm@rk\enspace}
	\parindent=\wd0
	\global\bibc@unt=0
	\def\bibitem##1{\global\advance\bibc@unt by 1
		\edef\lastref{\number\bibc@unt}
		{\toks0={##1}
		\message{bib[\lastref] \the\toks0,}%
		\immediate\write\@auxout
		{\noexpand\bibcite{\the\toks0}{\lastref}}}
		\def\next{##1}%
		\expandafter\ifx
		\csname bib@@\meaning\next\endcsname\lastref
		\else\global\bib@changestrue\fi%
		\bibcite{##1}{\lastref}
		\medbreak
		\item{\hfill\bibleftm@rk\lastref\bibrightm@rk}%
		}
	}
\def\endthebibliography{\egroup\par}
%
% THE NEXT MACRO MUST BE INCLUDED
% IN THE \BYE COMMAND. FOR INSTANCE:
%
    %\catcode`@=11
    \outer\def\bye{\@closeaux
    	\par\vfill\supereject\end}
    %\catcode`@=12
%
\def\@closeaux{\closeout\@auxout
	\ifref@changes\immediate\write16
	{Warning, changes in references.}\fi
	\ifpageref@changes\immediate\write16
	{Warning, changes in page references.}\fi
	\ifbib@changes\immediate\write16
	{Warning, changes in bibliography.}\fi
	\ifref@undefined\immediate\write16
	{Warning, references undefined.}\fi
	\ifpageref@undefined\immediate\write16
	{Warning, page references undefined.}\fi
	\ifbib@undefined\immediate\write16
	{Warning, citations undefined.}\fi}
%
% initialization of jobname.aux
%
\immediate\openin\@auxout=\jobname.aux
\ifeof\@auxout \immediate\write16
     {Creating file \jobname.aux}
\immediate\closein\@auxout
\immediate\openout\@auxout=\jobname.aux
\immediate\write\@auxout {\relax}%
\immediate\closeout\@auxout
\else\immediate\closein\@auxout\fi
%
% Let's read this file and open it out
%
\input\jobname.aux \par
\immediate\openout\@auxout=\jobname.aux
% this file will be closed by \bye.
%
% That's all, folks!
%

\def\cO{ {\cal O}}

\def\bR{{\bf R}}

\def\bx{{\bf x}}

\def\by{{\bf y}}

\def\bz{{\bf z}}
\def\bZ{{\bf Z}}

  \def\pro{\noindent {\bf{Proof.
}}}

\def\om{{\omega}}

\def\and{\quad\hbox{and }\quad}
\def \lrcorner
{\mathrel{\hbox{\vrule depth 0 pt height 0.4pt width 4pt
\vrule depth 0 pt height 5pt width 0.4pt
\kern 1pt }}}

\def\build#1_#2^#3{\mathrel{\mathop{\kern 0pt#1}\limits_{#2}^{#3}}}

%\newfam\gothfam \scriptscriptfont\gothfam=\fivegoth
%\textfont\gothfam=\tengoth \scriptfont\gothfam=\sevengoth
%\def\goth{\fam\gothfam\tengoth}

%\def\cqfd{\unskip\kern 6pt\penalty 500 \raise 0pt\hbox{\vrule\vbox
%to6pt{\hrule width 6pt \vfill\hrule}\vrule}\par}

\def\pro{\noindent {\bf Proof. }}

\def\smallsquare{\vbox{\hrule\hbox{\vrule height 1 ex\kern 1
ex\vrule}\hrule}}
\def\cqfd{\hfill \smallsquare\vskip 3mm}

\def\cC{{\cal C}}

\def\bibn@me{R\'ef\'erences bibliographiques}
%\input typpo
%
%\catcode`@=11
\def\bibliographym@rk{\bgroup}
%
% \bye est modifie pour la biblio et la table des matieres
%
\outer\def\bye{ 	\par\vfill\supereject\end}

%%%%%%%%%%%%%%%%%%%%%%%%%%%%%%%%%%%%%%%%%%

\null

\vskip 2mm

\centerline{\bf Approximation to points in  the plane by ${\rm SL(2,\bZ)}$-orbits}

\vskip 8mm

\centerline{Michel  L{\sevenrm AURENT} 
\& Arnaldo  N{\sevenrm OGUEIRA}
}

\vskip 11mm

\section{1. Introduction and results}
We view the real plane $\bR^2$ as a space of column vectors on which the 
 lattice  $\Gamma ={\rm SL(2,\bZ)}$ acts by left multiplication.
 Let $\bx =  \left(\matrix{ x_1\cr x_2\cr}\right)$ be a point in $ \bR^2$ with irrational slope $\xi = x_1/x_2$. 
The orbit $\Gamma \bx$ is then dense in $\bR^2$. The assertion follows from     J-S Dani's  density results \cite{Da}
for lattice orbits in homogeneous spaces, see also a more elementary proof  in   \cite{DaNo}.
The study  of lattice orbit distribution has been   the subject of numerous  works, in particular  \cite{MaWe},  
 \cite{No02} and \cite{No03} are concerned in counting the number of elements  $\gamma \bx$ belonging to various  sets under restriction on
  the size of $\gamma$, and  \cite{Gu} regards the  approximation  to radius with rational slope. 
Here we are concerned with the effective   approximation of a given point $\by\in \bR^2$ by points of  the form $\gamma \bx$, where $\gamma\in \Gamma$, in terms  of the size  of $\gamma$.

As a guide to  our results, let us   recall some classical  results of inhomogeneous approximation in $\bR$. 
 Minkowski Theorem asserts that for any irrational number $\xi$ and any real number $y$  not belonging to $\bZ \xi + \bZ$, 
 there exist infinitely many pairs of integers $(u,v)$, with $v \not= 0$, such that 
$$
| v \xi + u - y | \le {1\over 4 | v | }. \leqno{(1.1)}
$$
Our first goal is to obtain an analogous result 
  for the orbit $ \Gamma \left(\matrix{ \xi\cr 1\cr}\right)$ in $\bR^2$.
  Let  equip $\bR^2$ with   the supremum  norm $|\bx |= \max(|x_1|,|x_2|)$, and for any matrix $\gamma$,  denote as well by $|\gamma |$
 the maximum of the absolute values of the entries of $\gamma$.  Notice that any choice of  norm on the algebra  of matrices 
 ${\rm M}_2(\bR)$  would lead to the same exponents with possibly different constants. 
 We distinguish three cases, according as the target point $\by$ coincides with the origin ${\bf 0}=  \left(\matrix{0\cr 0\cr}\right)$, or it lies on a radius whose slope
 is either rational or irrational.
 
 \medskip
 \noindent
{\bf Theorem 1.} {\sl  Let  $\bx$ be a  point in $\bR^2$ with irrational slope. 

\noindent
{\rm (i)} There exist infinitely many  matrices $\gamma \in \Gamma$ such that
$$
|  \gamma \bx |  \le { | \bx |   \over | \gamma |} . \leqno{(1.2)}
$$
{\rm (ii)} Let   $\by =  \left(\matrix{ y_1\cr y_2\cr}\right)$ be a point $ \in \bR^2\setminus\{{\bf 0}\}$.  Assume  that 
either the slope  $ y_1/y_2$  is  a rational number  $a/b$, or that $y_2=0$  in which case we put
$a=1$ and $b=0$. Then,  there exist infinitely many  matrices $\gamma \in \Gamma$ such that
$$
| \gamma \bx - \by | \le {c \over | \gamma |^{1/2}} \quad {\rm with} \quad  c =   2\sqrt{3}\max( | a |, | b |) | \bx|^{1/2} | \by |^{1/2}.\leqno{(1.3)}
$$
{\rm (iii)}  When the slope $y_1/y_2$ of the point $\by$ is irrational, there exist  infinitely many  matrices $\gamma \in \Gamma$ satisfying
 $$
 | \gamma \bx - \by | \le {c' \over | \gamma |^{1/3}} \quad {\rm with} \quad  c' =   7\sqrt{5} | \bx|^{1/3} | \by |^{2/3}.\leqno{(1.4)}
$$
}

The exponents 1 and $1/2$ occurring respectively in (1.2) and (1.3) are best possible. We are also interested in {\it uniform}
versions of Theorem 1,  in the sense of \cite{BuLa}. We  first state the uniform version of Minkowski Theorem. To this purpose,
we need  the standard   notion   of {\it  irrationality measure} of an irrational number.  

\proclaim
Definition. For any irrational real number $\alpha$, we denote by   $\om(\alpha)$ the supremum of the numbers $\om$ such that the inequation 
$$
\vert v \alpha   +u   \vert \le \vert v \vert^{-\om}
$$
has infinitely many integer solutions  $(v,u)$. 

Then, for any  real number $\mu < 1/\om(\xi)$ and  any positive real number $T$ sufficiently large in terms of $\mu$, there
exists integers $u,v$ such that 
$$
\max( | u |,  | v|)  \le T \and | v \xi + u - y | \le   T^{-\mu}. \leqno{(1.5)}
$$
See for instance  the main theorem of \cite{BuLa},  as well as the comments  explaining  the 
link with the claims (1.1) and (1.5).    
 More information and results can be found in \cite{Bu,BuLa,Cas},
 including metrical theory and higher dimensional generalizations. 
 
 In view of the above results,  let us define   two  exponents  
 measuring respectively the usual and  the uniform approximation to a  point $\by$ by elements of the  orbit $\Gamma \bx$.
 We follow the notational conventions of \cite{BuLa}.

\proclaim
Definition. Let  $\bx$ and  $\by$ be two points in   $\bR^2$. We  denote by $\mu(\bx, \by)$ the supremum of the real numbers  $\mu$
for which there  exist infinitely many matrices  $\gamma \in \Gamma$ satisfying the inequality 
$$
| \gamma \bx - \by | \le | \gamma |^{-\mu}.
$$
We  denote by $\hat\mu(\bx, \by)$ the supremum of the exponents $\mu$ such that for any sufficiently large positive real number $T$, 
there  exists  a  matrix   $\gamma \in \Gamma$ satisfying 
$$
 | \gamma | \le T \and | \gamma \bx - \by | \le T^{-\mu}.
$$

Clearly $ \mu(\bx,\by) \ge \hat\mu(\bx,\by)\ge 0$,  unless $\by$ belongs to the orbit $\Gamma\bx$  in which case $\hat\mu(\bx,\by)= + \infty$. We can now state the
 
  \medskip
 \noindent
{\bf Theorem 2.} {\sl  Let  $\bx$ be a  point in $\bR^2$ with irrational slope $\xi$. 

\noindent
{ \rm (i)} We have  
$$
\mu(\bx,{\bf 0}) = 1 \and \hat\mu(\bx,{\bf 0}) ={1\over \om(\xi)}. \leqno{(1.6)}
$$
{Ê\rm (ii)}  Let   $\by =  \left(\matrix{ y_1\cr y_2\cr}\right)$ be a point $ \in \bR^2\setminus\{{\bf 0}\}$. Assume  that either the  slope $ y = y_1/y_2$
is  rational or that $y_2=0$. Then,   we have  the equalities
$$
\mu(\bx,\by) =  {\om(\xi)\over \om(\xi) + 1}  \ge {1\over 2} \and \hat\mu(\bx,\by) = {1 \over  \om(\xi)+ 1} . \leqno{(1.7)}
$$
{\rm (iii)} When the slope $y$ of the point $\by$ is an  irrational number, then the following  lower bounds hold
$$
\mu(\bx,\by) \ge {1\over 3} \and \hat\mu(\bx,\by) \ge { \om(y) +1\over 2(2 \om(y)+1) \om(\xi)} \ge {1 \over 4 \om(\xi)}. \leqno{(1.8)}
$$
}

\bigskip 
If $\xi$ is  a Liouville number, meaning that $\om(\xi) = +\infty$, the equalities (1.7) obviously read $ \mu(\bx,\by) =1$ and $\hat\mu(\bx,\by) =0$.
When the  slope $y$ is rational, an explicit lower bound for the distance between $\gamma \bx$ and $\by$ will be given in Theorem 4 of Section 8,
which brings further information in terms  of the convergents of $\xi$. 

Note that Maucourant and Weiss \cite{MaWe} have recently obtained the weaker lower bounds 
$$
\mu(\bx,\by) \ge {1\over  144}  \and
  \hat\mu(\bx,\by) \ge {1 \over 72 (\om(\xi)+1)}, 
  $$
  as a consequence of  effective equidistribution estimates for unipotent trajectories in $\Gamma\backslash{\rm SL(2,\bR)}$ (use  Corollary 1.9 
 in  \cite{MaWe}  and substitute $\delta_0 = 1/48$,  which is an admissible value as mentioned in Remark 1.6).
 In  another  related work \cite{Gu},  Guilloux observes the existence of gaps around rational directions
 in the repartition of the {\it cloud}  of points $\{\gamma\bx \, ;  \gamma \in \gamma , | \gamma | \le T\}$ for large $T$.
In our setting,  he  proves the upper  bound $\hat\mu(\bx,\by) \le 1$ for any point $\by$ with rational slope.

  We now discuss  upper bounds for our exponents $\mu(\bx,\by)$ and $\hat\mu(\bx,\by)$.
   Applying  Proposition 8 of \cite{BuLa} to the two  inequalities  of the form  (1.5)
   determined   by the two   coordinates   of $\gamma \bx -\by$, 
   we obtain  the bound  
  $
  \hat\mu(\bx,\by) \le \om(\xi)
  $
 for any point $\by$ which   does not belong to the orbit $\Gamma\bx$.  
  Moreover, the stronger upper bound
 $$
 \hat\mu(\bx,\by) \le {1\over \om(\xi)} \le \om(\xi)
 $$
 holds for almost all \footnote{(*)}{\sevenrm Throughout the paper, the expression `almost all'
 always refers to Lebesgue measure in the ambient space.}
 points $\by$, since the main theorem   of \cite{BuLa} tells us that the exponent $\mu $ 
 in (1.5) cannot be larger than $1/\om(\xi)$ for almost all real number $y$. 
  As for the exponent $\mu(\bx,\by)$, it  may be arbitrarily large when  $\by$ is a point of {\it Liouville} type, meaning that $\by$ is the limit of 
  a fast converging  sequence $(\gamma_n\bx)_{n\ge 1}$ of points of the orbit. However,  $\mu(\bx,\by)$ is  bounded  almost
  everywhere. Projecting as above on both coordinates,  the main theorem of \cite{BuLa} shows that the upper bound 
  $
  \mu(\bx, \by) \le 1
  $
  holds for almost all points $\by$. Here is a stronger statement.  
  
 \proclaim
 Theorem 3.
 Let  $\bx$ be a point in  $\bR^2$ with irrational slope  and let $y$ be an irrational  number having   irrationality measure 
 $\om(y)=1$.
 Then,  the upper bound
 $$
 \mu(\bx,  \by) \le {1 \over 2}
 $$
holds for almost all  points $\by $ of  the line $\bR  \left(\matrix{ y  \cr 1 \cr}\right)$.

 It follows from  theorems 2 and 3 that, $\bx$ being fixed,  we have  the estimate
 $$
 {1\over 3} \le \mu(\bx, \by)  \le { 1\over 2}
 $$
   for  almost all  points $\by \in \bR^2 $, since the assumption $\om(y) = 1$ occurring in Theorem 3 is valid  for almost all real numbers  $y$. 
 Moreover the maximal  value $1/2$   is reached
 for  any  point $\by\not= {\bf 0} $ lying on a radius with rational slope when the slope $\xi$ of $\bx$  has  irrationality measure
 $\om(\xi) =1$.  We adress the problem of finding the generic value, 
if it does exist,   of the exponents
$\mu(\bx,\by)$ and $\hat\mu(\bx,\by)$ on $\bR^2\times \bR^2$.
Heuristic (but optimistic) equidistribution arguments suggest that we should have
$$
\mu(\bx,\by) = \hat\mu(\bx,\by) ={1\over 2}
$$
for almost all pairs of points $(\bx,\by)$.

 Let us detail the content of the paper.   In Section 2, we associate to an irrational  number $\xi$  
 a  sequence of  matrices in $\Gamma$, called {\it convergent matrices},  which  send any  point  $\bx$ with  slope $\xi$ 
 towards the origin. As first application, the  easy case   $\by ={\bf 0}$  is   investigated in Section 3.
 In Section 4, we expand tools for constructing approximants  to a point $\by$ by elements $\gamma \bx$ of the orbit. 
 Our approach is explicit. 
  We write down  $\gamma$ 
 as a product of three factors $NGM$. The matrix $M$ is some  convergent matrix associated to the slope  $\xi$ of $\bx$, while the 
 matrix $N$ is essentially the inverse of a convergent matrix associated to  the slope of the target point $\by$. 
  As for the factor $G$, whose choice is not uniquely determined, we use some  suitable unipotent matrix. 
 From a dynamical point of view, the way for going from $\bx$ to
 $\by$ splits  into three different stages. First, we push  down  $\bx$ close to the origin.   
 Next, we move on an horizontal line (any fixed rational
 direction should be convenient), and finally we point   up to $\by$ thank to the third factor $N$.
  We apply the method in Sections 5 and 6, thus obtaining  various  lower bounds for $\mu(\bx, \by)$ and $\hat\mu(\bx,\by)$ 
 depending on whether the slope of the point $\by$ is rational or not.  On the other hand,  
  we obtain  upper bounds for these exponents in Sections 7 and 8.  Conversely, 
 a decomposition of the  form $\gamma = NGM$, with a factor $G$ of small norm,  is in fact necessary; it  implies  upper bounds valid for almost all points
 $\by \in \bR^2$,  including all points $\by$  with rational slope. 
  In the latter case, it turns out that the upper and lower bounds thus obtained coincide; hence we get exact values
 for $\mu(\bx,\by)$ and $\hat\mu(\bx,\by)$.
 The final Section 9  deals with additional constraints of signs.
 
 It would   be interesting to extend our decomposition method to other lattices $\Gamma$ in ${\rm SL(2,\bR)}$. Observe that
 the rational slopes, namely the  cusps of the Fuchsian group $ {\rm PSL(2,\bZ)}$, play  a prominent role in our approach.

 We write $A \ll B$ when  there exists a positive constant $c$ such that $A \le c B$ for all values of the 
 parameters   under consideration (usually the indices $j$ and $k$).
The coefficient $c$ may possibly depend upon the points $\bx$ and $\by$. 
 As usual, the notation  $A \asymp B$ means that $A \ll B$ and $A \gg B$.

\section
{2. Convergent matrices  }
 Let  $\xi$ be an irrational number and let $(p_k/q_k)_{k\ge 0}$ be the sequence of convergents of  $\xi$. We set  
  $
  \epsilon_k = q_k\xi - p_k.
  $
  The theory of continued fractions tells us that the sign of $\epsilon_k$ is alternatively positive or  negative according to whether   $k$ 
  is even or odd, and that the estimate 
  $$
  {1\over 2 q_{k+1} } \le | \epsilon_k | \le {1\over q_{k+1}}  \leqno{(2.1)}
  $$
  holds for $k\ge 0$. For later use,  note  as a consequence of (2.1) that,  when $\om(\xi)$ is finite, we have the upper bound
  $q_{k+1} \le q_k^\om$ for any real number $\om > \om(\xi)$ provided $k$ is large enough, while if
  $\om < \om(\xi)$, the lower bound $q_{k+1} \ge q_k^\om$ holds for infinitely many $k$. 
  
For any positive integer $k$, we set
$$
M_k =  \left(\matrix{ q_k & -p_k\cr -q_{k-1}& p_{k-1}\cr}\right) \quad {\rm or} \quad  
M_k =  \left(\matrix{ q_{k} & -p_{k}\cr q_{k-1}& - p_{k-1}\cr}\right) , 
 $$
respectively when     $k$ is even or odd. In both cases the  matrix $M_k$  belongs  to $\Gamma$ and has  norm $| M_{k} |  = \max(q_{k}, | p_{k} | )$.
 Let 
$\bx =  \left(\matrix{ x_1\cr x_2\cr}\right)$ be a point with slope $\xi= x_1/x_2$. Then,    we have
$$
M_k\bx =  x_2\left(\matrix{\epsilon_{k} \cr (-1)^{k-1}  \epsilon_{k-1}  \cr}\right)=  x_2\left(\matrix{\epsilon_{k} \cr |  \epsilon_{k-1} |   \cr}\right),
 $$
noting that the second coordinate $(-1)^{k-1}\epsilon_{k-1}$ is  always positive
and thus equals  $| \epsilon_{k-1} |$. 

The matrices   $M_{k}$
will be    called  {\it  convergent  matrices} of $\xi$. The name is justified by the fact that  the  columns  of the  inverse matrix 
 $$
 M_{k}^{-1} =  \left(\matrix{ p_{k-1} &  p_{k}\cr  q_{k-1}& q_{k}\cr}\right)
 \quad {\rm or} \quad  M_{k}^{-1} =  \left(\matrix{ -p_{k-1} &  p_{k}\cr  -q_{k-1}& q_{k}\cr}\right)
$$
 give,  up to a sign, the numerator and the denominator of   two consecutive convergents  of  $\xi$.

 \section
 {3. Approximation to the origin}
  We first consider  the  easier case where the  target point $\by$ equals  the origin ${\bf 0}$, and 
  prove   in this section  the claims (1.2) and (1.6). We assume without loss of generality that $\bx =  \left(\matrix{ \xi \cr 1 \cr}\right)$.  
 
   \proclaim
  Lemma 1. Let $k$ be a positive integer and let $\gamma \in \Gamma$ with norm $ | \gamma | \le q_{k+1}/2$. 
  Then, we have the lower bound 
  $$
  | \gamma \bx | \ge {1\over 2 q_k}.
  $$
 
 \pro We argue by contradiction. On   the contrary, suppose  that   
 $
  | \gamma \bx | <  {1/( 2 q_k)}.
  $
Put $\gamma = \left(\matrix{ v_1 & u_1\cr v_2& u_2\cr}\right)$ and  $G = \gamma M_k^{-1}$.  Assume first that $k$ is even. 
We find the formula
$$
\eqalign{
G = &  \left(\matrix{ v_1 & u_1\cr v_2& u_2\cr}\right) \left(\matrix{ q_{k} & -p_{k}\cr -q_{k-1}& p_{k-1}\cr}\right)^{-1}
 =   \left(\matrix{  p_{k-1}v_1 + q_{k-1} u_1&  p_{k}v_1  + q_{k} u_1\cr   p_{k-1}v_2 + q_{k-1} u_2&   p_{k}v_2 + q_{k} u_2\cr}\right)
 \cr
 = & \left(\matrix{- v_1( q_{k-1}\xi - p_{k-1})   + q_{k-1}( v_1\xi +u_1)  & - v_1( q_{k}\xi - p_{k})   + q_{k}( v_1\xi +u_1) 
  \cr - v_2( q_{k-1}\xi - p_{k-1})   + q_{k-1}( v_2\xi +u_2) & - v_2( q_{k}\xi - p_{k})   + q_{k}( v_2\xi +u_2) \cr}\right) .
 \cr}
 $$
 Bounding  from above the norm of the second column of the above matrix gives   
$$
\max \Big( | - v_1( q_{k}\xi - p_{k})   + q_{k}( v_1\xi +u_1) |  , |  - v_2( q_{k}\xi - p_{k})   + q_{k}( v_2\xi +u_2) |\Big) \le
{ | \gamma | \over q_{k+1}} + q_k | \gamma \bx | <1.
$$
Since $G$ has integer entries, it follows that the second column of $G$ equals  ${\bf 0}$. 
The case $k$ odd leads to the same conclusion. Contradiction   with $\det G =1$. 
\cqfd

\medskip
For any $\gamma \in \Gamma$ of norm $| \gamma | >  q_1/2 $, let $k$ be the integer defined  by the estimate
$$
{q_k \over 2 } <  | \gamma  | \le { q_{k+1}\over 2}.
$$
It follows from Lemma 1 that
$$
| \gamma \bx | \ge {1 \over 2 q_k} \ge {1 \over 4 | \gamma |}.
$$
Therefore $ \mu(\bx , {\bf 0}) \le 1$. On the other hand, we have  that 
 $$
 | M_k | = \max( | p_k | ,   q_k)  \and  | M_k  \bx | = \max( | \epsilon_{k} | , | \epsilon_{k-1} | )= | \epsilon_{k-1} |  \le  { 1\over q_k }, 
 $$
 by (2.1). Observe that $p_k = q_k \xi - \epsilon_k$ has absolute value $\le | \xi |  q_k$ if  $\epsilon_k$ and $\xi$ have the same sign.
 Hence (1.2) holds for $\gamma = M_k$ when $k$ is either odd or even. 
 
 It obviously  follows from (1.2) that  $\mu(\bx, {\bf 0}) =1$, thus proving the first assertion of (1.6). 
 The proof of the equality $\hat\mu(\bx,{\bf 0}) = 1/\om(\xi)$ is similar. For any real number  $\om < \om(\xi)$, 
 there exist  infinitely many $k$ such that
$ q_{k+1} \ge q_k^\om$. Put $T = q_{k+1}/2$. For all $\gamma \in \Gamma$ with norm $| \gamma | \le T$, Lemma 1 gives the lower bound
$$
 | \gamma \bx | \ge {1\over 2 q_k } \ge {1\over 2 (2 T)^{1/\om}}.
 $$
  Therefore $ \hat\mu(\bx, {\bf 0}) \le 1/\om$,  
 and letting $\om$ tend to $\om(\xi)$, we obtain  the upper  bound $\hat\mu(\bx , {\bf 0}) \le 1/\om(\xi)$.
On the other hand,  the choice of the matrix $\gamma = M_k$ for $ | M_k | \le T  <  | M_{k+1} |$ shows  that $\hat\mu(\bx , {\bf 0}) \ge 1/\om(\xi)$.
Hence the equality $\hat\mu(\bx , {\bf 0}) = 1/\om(\xi)$ holds.

  \section
{4. Construction of approximants} 
The group $\Gamma$ is generated by the two matrices
$$
J =   \left(\matrix{ 0 & -1\cr 1& 0\cr}\right)
\and 
U=   \left(\matrix{ 1 & 1\cr 0& 1\cr}\right)
$$
satisfying the relations $J^2 = -{\rm Id}$ and $(JU)^3 = -{\rm Id}$. Observe  that the matrix $J$ acts on $\bR^2$ as a rotation of a right angle,
while the unipotent matrix $U$ leaves invariant any horizontal line $ \left\{\left(\matrix{ z \cr \epsilon\cr}\right) ; z\in \bR\right\}$
  and acts on this line as a translation 
of step $ \epsilon $.

 From now on, we   assume that the target 
 point $\by$ differs from ${\bf 0}$. Note that   $| J\bz | = | \bz | $ for all $\bz \in \bR^2$. 
Replacing possibly  $\bx$ by  $J \bx $  or  $\by$ by  $ J \by$,  we shall assume throughout  the paper that 
$$
| \bx | = | x_2 | \and | \by | =   | y_2 |,
$$
 so that the slopes $\xi = x_1/x_2$ and $y = y_1/y_2$  of the points $\bx$ and $\by$ satisfy 
 $$
 0 <  | \xi | <  1 \and    | y | \le 1.
  $$
  
We consider  matrices of the form
 $
\gamma =  N U^\ell M_{k},
 $
 where $\ell $ is an integer and  $N$ is a  matrix in $\Gamma$, which will be specified later.

 \proclaim
 Lemma 2. Let  $k$ be a positive integer,   $\ell$ be an integer,  and let  $N=  \left(\matrix{ t & t'\cr s & s' \cr}\right)$ belong to  $\Gamma$. 
 Put  $\gamma = NU^{\ell} M_k \in \Gamma$.  Then
  $$
\left\vert \ell q_{k-1} + (-1)^{k-1}q_k \right\vert  | s |   - | s' | q_{k-1} \le 
 | \gamma | \le  | \ell | | N | q_{k-1} +  2 | N |  q_k.
 $$
 
 \pro  Since  $| \xi |  <1$, we have  $ |p_k |\le q_k$  for all  $ k\ge 0$. 
 When   $k$ is even, we have  
$$
\eqalign{
\gamma  = & \left(\matrix{ t & t'\cr s & s' \cr}\right) \left(\matrix{1 & \ell \cr 0&1\cr}\right) \left(\matrix{ q_k & -p_k\cr -q_{k-1}& p_{k-1}\cr}\right)
\cr  
= & \left(\matrix{ - \ell t q_{k-1} +  t q_k - t'q_{k-1} &  \ell t p_{k-1} -  t p_k + t'p_{k-1} 
\cr  - \ell s q_{k-1} + s q_k - s'q_{k-1}  & \ell s p_{k-1} - s p_k + s'p_{k-1} \cr}\right) .
\cr}
$$
 When $k$ is odd, we find 
 $$
\eqalign{
\gamma  = & \left(\matrix{ t & t'\cr s & s' \cr}\right) \left(\matrix{1 & \ell \cr 0&1\cr}\right) \left(\matrix{ q_k & -p_k\cr q_{k-1}& -p_{k-1}\cr}\right)
\cr  
= &  \left(\matrix{  \ell t q_{k-1} +  t q_k + t'q_{k-1} &  -\ell t p_{k-1} -  t p_k - t'p_{k-1} 
\cr   \ell s q_{k-1} + s q_k + s'q_{k-1}  & -\ell s p_{k-1} - s p_k - s'p_{k-1} \cr}\right) .
\cr}
$$
The required upper bound obviously holds in both cases. For the lower bound,  look at the lower left entry of $\gamma$. \cqfd

 \proclaim
 Lemma 3.   Let  $k$ be a positive integer,   $\ell$ be an integer,  let  $N = \left(\matrix{ t & t' \cr s & s' \cr}\right)$ belong to  $\Gamma$
and let $y$ be any  real number.
 Put 
 $$
\gamma =  N U^{\ell} M_k =   \left(\matrix{ v_1 & u_1\cr v_2& u_2\cr}\right), \quad \delta = | sy-t| \and  \delta' = | s'y -t' |.
$$
 Then, we have the upper bound
 $$
 | v_1 \xi + u_1 -y(v_2 \xi +u_2) | \le { \delta | \ell | \over q_k} + {\delta \over q_{k+1}} + { \delta' \over q_k}.
 $$

\pro  It is a simple matter of bilinearity. We have the formula
$$
\eqalign
{
y(v_2\xi + u_2)-v_1\xi-u_1  = &\left(\matrix{ - 1 &  y \cr}\right) \gamma
 \left(\matrix{ \xi  \cr 1 \cr}\right)
\cr
=&  \left(\matrix{ - 1 &  y \cr}\right)
 \left(\matrix{ t & t' \cr s & s'\cr}\right) \left(\matrix{1 & \ell \cr 0 & 1\cr}\right) M_k \left(\matrix{ \xi  \cr 1 \cr}\right)
\cr 
= & \left(\matrix{ s y -t &  s'y-t'  \cr}\right)  \left(\matrix{1 & \ell \cr 0 & 1\cr}\right)\left(\matrix{ \epsilon_k  \cr | \epsilon_{k-1}|  \cr}\right)
\cr
= & (sy-t)(\epsilon_k +    \ell | \epsilon_{k-1} | ) +  (s'y-t')| \epsilon_{k-1}|.
\cr}
$$
 Now the upper bound immediately follows from the estimate (2.1).
 \cqfd

 We shall use Lemma 3 in the following way.   Put
 $$
 \left(\matrix{ \Lambda_1 \cr  \Lambda_2 \cr}\right)= \gamma \bx - \by =  \left(\matrix{  x_2(v_1\xi +u_1) -y_1 \cr   x_2(v_2 \xi +u_2) -y_2 \cr}\right) 
 $$
 and  let $y= y_1/y_2$ be the slope of the  point $\by$, so that  
 $$
 \Lambda_1 - y \Lambda_2 =  x_2\Big(v_1\xi +u_1 -y(v_2\xi + u_2)\Big).
 $$
 Now,  Lemma  3 provides us with a  fine  upper bound for $ |  \Lambda_1 - y \Lambda_2 | $,  
 as far as the quantities $\delta$ and $\delta'$ are small. Therefore to bound from above $| \gamma \bx - \by |$, 
   it suffices to bound one of its coordinates,
  say $\Lambda_2$.     We have the expression
  $$
 \Lambda_2 = x_2 \Big( s \epsilon_k +  (s\ell+s')| \epsilon_{k-1}| \Big) -y_2 =  x_2 s | \epsilon_{k-1} | ( \ell - \rho) , \leqno{(4.1)}
 $$
 where 
 $$
 \rho = {y_2\over x_2 s | \epsilon_{k-1}| } - {\epsilon_k\over | \epsilon_{k-1}| } - {s'\over s}. \leqno{(4.2)}
 $$

\section
{4.1.  Irrational slopes }
 We assume here that the slope  $y=y_1/y_2$ is an irrational number and  apply the key lemmas 2 and 3 for  constructing 
matrices   $\gamma$  in $\Gamma$ such that $\gamma \bx $ is close to $\by$.
 
 Denote by   $(t_j/s_j)_{j\ge 0}$  the sequence of convergents of   $y$, and put 
$$
 N_j =  \left(\matrix{ t_j & t'_j\cr s_j& s'_j\cr}\right) ,  \quad {\rm where} \quad  s'_j =(-1)^{j-1}s_{j-1} \and  t'_j =  (-1)^{j-1}t_{j-1}, 
$$
for any $j\ge 1$.  Observe that $J N_j^{-1}$ coincides with the convergent matrix
$M_j $ associated to the irrational number $y$ as in  Section 2. Hence  $N_j$ belongs to $\Gamma$. 
 
 \proclaim
 Lemma 4. Let $j$ and $k$ be positive integers. There exists a matrix $\gamma\in \Gamma$, of the form $N_jU^{\ell}M_k$
 for some integer $\ell$, such that
 $$
\left\vert {  | y_2 | \over | x_2 |} q_{k-1}q_k -  s_j q_k \right\vert  -4 s_jq_{k-1} \le 
 | \gamma | \le { 2 | y_2 | \over | x_2 |} q_{k-1}q_k + 4 s_j q_k  \leqno{(4.3)}
 $$
and
$$ 
  |  \gamma \bx -\by | \le { 2 | y_2 | \over s_j s_{j+1}} +  { 5 | x_2 | s_j \over q_k}. \leqno{(4.4)}
 $$
 
  \pro
 Since  $ | y | <   1$, we have  $ | t_j | \le s_j$ and  $| t'_j |  \le | s'_j |  < s_j$.
 The matrix  $N_j$ has thus  norm 
 $
 | N_j |   = s_j.
 $
 The theory of continued fractions gives the upper bounds
 $$
 \delta = | s_j y - t_j | \le s_{j+1}^{-1} \and \delta' = |  s'_j y - t'_j | = | s_{j-1}y - t_{j-1} | \le s_j^{-1}. \leqno{(4.5)}
 $$
 Recall the definition of  $\rho$  given in  (4.2), and  substitute  $s_j$ to $s$ and $s'_j$ to  $s'$.
 Bounding $| \epsilon_k/\epsilon_{k-1} | \le 1$,  $s_{j-1}/s_j \le 1$, and 
$q_k \le | \epsilon_{k-1} |^{-1} \le 2 q_k$ by (2.1), we find
$$
 {  | y_2 |  q_k \over | x_2  | s_j } -2   \le  |  \rho | \le { 2 | y_2 |  q_k \over | x_2  | s_j } + 2. 
$$
Define  $\ell$ as being  the unique integer such that
$$
 | \ell - \rho | <  1 \and | \ell | \le | \rho | .  
$$
We  set  
$$
\gamma = N_jU^{\ell}M_k \and   \left(\matrix{ \Lambda_1\cr \Lambda_2 \cr}\right) = \gamma \bx -\by.
$$
Therefore 
$$
 {  | y_2 |  q_k \over | x_2  | s_j } - 3  \le  | \ell  | \le  { 2 | y_2 |  q_k \over | x_2  | s_j } + 2,\leqno{(4.6)}
$$
and it follows from    (4.1) that 
$$
 |  \Lambda_2  |   
 = | x_2 |  s_j | \epsilon_{k-1}| |  \ell - \rho| \le {  | x_2 | s_j \over q_k } . 
$$
Now, we apply  Lemma 3 to bound $\Lambda_1 - y \Lambda_2$. Using (4.5) and (4.6), we find
$$
| \Lambda_1 - y  \Lambda_2 | \le  | x_2 | \left( { | \ell | \over s_{j+1}q_k } + { 1 \over s_{j+1} q_{k+1} } + { 1 \over s_j q_k } \right) 
\le | x_2 | \left( { 2 | y_2  | \over | x_2 | s_j s_{j+1} } + { 4 \over s_j  q_k}   \right) .
$$ 
Since $| y| <  1$, summing   the two above upper bounds  gives
$$
| \Lambda_1 | \le   |\Lambda_2 | + | \Lambda_1 - y  \Lambda_2 | 
 \le | x_2 | \left( { 2 | y_2  | \over | x_2 | s_j s_{j+1} } +  {  5 s_j \over q_k } \right).
$$
We have obtained the upper bound 
$$
| \gamma \bx -\by | = \max ( | \Lambda_1 |, | \Lambda_2 | )   \le   { 2 | y_2 | \over s_j s_{j+1}} +  { 5 | x_2 | s_j \over q_k}
$$
claimed in (4.4). On the other hand, Lemma 2 combined with     (4.6)  gives  the estimate  of  norm (4.3). 
\cqfd

\section{ 4.2 Rational slopes}
We  consider here  a target point  $\by$ with rational slope $y$.  Writing the rational $y = a/b$ 
in reduced form,  the integers  $a$ and $b$ are  coprime and we have     $| a | \le b$, 
since we have assumed that $ |  y  | \le 1$.

\proclaim
Lemma 5. For any sufficiently large integer $k$, there exists a matrix $\gamma \in \Gamma$ such that
$$
 { | y_2 | \over 2 | x_2 |} q_{k-1}q_k  \le   | \gamma | \le   { 3| y_2 | \over | x_2 |} q_{k-1}q_k 
\and
| \gamma \bx -\by |  \le { 2 b |x_2 |\over  q_k }.
$$

\pro
We now use  as best rational approximation to 
 $y$ the number  $y= a/b$ itself.
 
Let us complete the primitive point $ \left(\matrix{ a \cr b\cr}\right)$
into  an  unimodular matrix 
$
N=  \left(\matrix{ a & a' \cr  b & b' \cr}\right),
$
with norm $| N |  = b$. The matrix $N$ is thus fixed, independently of $k$, and we have 
$$
\delta  = | by -a | = 0 
\and 
\delta'  = | b'y -a' | = {1\over b}. \leqno{(4.7)}
$$
We use lemmas 2 and 3 with  this choice of matrix $N$. 
Recall the definition of $\rho$ given in (4.2),
with $s$ and $s'$ respectively replaced by $b$ and $b'$.   As previously, define  $\ell$ as the unique integer verifying   
$| \ell | \le | \rho |$   and $ | \ell - \rho | < 1$.
We have the  estimate
$$ 
\left({  | y_2 |  \over b  | x_2 |   }\right) q_k -  3   \le | \ell |  \le  \left({2  | y_2 |  \over b | x_2 |   }\right) q_k + 2,  \leqno{(4.8)}
$$
 and 
$$
 |  \Lambda_2  |  
 = | x_2 |  b | \epsilon_{k-1}| |  \ell - \rho| \le { | x_2 | b \over  q_k }  .  \leqno{(4.9)}
$$
Substituting the values of $\delta$ and $\delta'$ given by (4.7), Lemma 3  now gives 
$$
| \Lambda_1 - y  \Lambda_2 | \le {| x_2 | \over b q_k}. \leqno{(4.10)}
$$
We deduce from (4.9), (4.10) and   the triangle inequality that 
$$
| \gamma \bx -\by |  \le   { 2 b |x_2 | \over  q_k }, 
$$
as claimed. Finally, taking (4.8) into account, Lemma 2 gives 
 $$
| \gamma | \le   |\ell | b  q_{k-1} + 2 b q_k  \le 2 { | y_2 | \over | x_2 |} q_{k-1}q_k + 2 b q_{k-1} + 2 b q_k 
\le    3{ | y_2 | \over | x_2 |} q_{k-1}q_k 
$$
and
$$
| \gamma | \ge   |\ell | b  q_{k-1} -  2 b q_k  \ge  { | y_2 | \over | x_2 |} q_{k-1}q_k - 5  b q_{k}  
\ge    { | y_2 | \over 2  | x_2 |} q_{k-1}q_k ,
$$
for large $k$. \cqfd

\section{5. Proof of  Theorem 1}
We  apply lemmas 4 and 5 in order to prove respectively the claims  (1.3) and (1.4).  
We first deal  with an  irrational slope $y$ and  prove (1.4) in the sections 5.1 and 5.2 below.
The argument  splits into two parts depending on whether   
 the value of the irrationality measure $\om(\xi)$ is smaller than 3 or greater than 2. 

\section
{5.1.  The case  $\om(\xi)  < 3 $}
 Let us  define infinitely many   pairs of integers  $j$ and  $k$ in the following way.   Let $j_0$ be an arbitrarily large integer.
We determine    $k$ by the estimate  
$$
\left( {| y_2 |  q_{k-1}\over | x_2 |  }\right)^{1/3} < s_{j_0} \le \left( {| y_2 |  q_{k}\over | x_2 |  }\right)^{1/3}.
$$
Let  $j$ be the largest integer such that  $s_j$ belongs to the above interval. We thus have the inequalities
$$
\left( {| y_2 |  q_{k-1}\over | x_2 |  }\right)^{1/3} < s_{j} \le \left( {| y_2 |  q_{k}\over | x_2 |  }\right)^{1/3}< s_{j+1}. \leqno{(5.1)}
$$

We use Lemma 4 for any  pair $j$ and $k$ verifying  (5.1). 
It  provides  us with a matrix $\gamma $ satisfying (4.3) and (4.4). Combining (4.4) and (5.1), we find  the upper bound
$$
| \gamma \bx -\by |   \le    | y_2 |^{1/3}  | x_2 |^{2/3} \left( {2 \over q_{k-1}^{1/3}q_k^{1/3}}
+ {5 \over q_k^{2/3}}\right) \le { 7 | y_2 |^{1/3}  | x_2 |^{2/3}\over (q_{k-1}q_k)^{1/3}} . \leqno{(5.2)}
$$
Observe now that for any  real number  $\om$ satisfying  $\om(\xi)< \om <3$, we have 
$
q_{k-1} \ge q_k^{1/\om} 
$
for all  $k$ sufficiently large. Since  $ s_j \ll q_k^{1/3}$, the second term  $4s_j q_k$ occurring on  the right hand side of  (4.3) is 
much smaller than  the first one, as  $k$  tends to infinity. Thus, for any sufficiently large  $k$, we have the  norm bound 
$$
| \gamma |  \le  3 { | y_2|\over | x_2|}q_{k-1}q_k  . \leqno{(5.3)}
$$
  Combining then  (5.2) and  (5.3), we obtain 
$$
| \gamma \bx -\by | \le 7 \root 3\of 3 | x_2|^{1/3}| y_2|^{2/3} |\gamma|^{-1/3}  \le   c' | \gamma |^{-1/3}.
$$

The  upper bound (1.4) is therefore established. It remains to show that our construction produces infinitely many solutions of (1.4).
To that purpose, it suffices to bound from below the norm of $\gamma$. The estimate  (4.3) in Lemma 4  gives indeed 
$$
| \gamma | \asymp   { | y_2 | \over |  x_2 |} q_{k-1}q_k.
$$

\section
{5.2.  The case  $\om(\xi) >2 $}
 Let us fix a real number  $\om $ satisfying   $2 < \om < \om(\xi)$. There exist infinitely many   $k$  such that  
$
q_{k-1}^\om  \le q_k . 
$
For any such  integer  $k$, let $j$ be the integer defined  by the inequality
 $$
 s_{j} \le \left( {| y_2 |  q_{k}\over | x_2 |  }\right)^{1/2}< s_{j+1}. \leqno{(5.4)}
 $$
% and let $i$ be the largest integer  such that 
% $$
%  s_{j-1} +i s_j  \le  \left( {| y_2 |  q_{k}\over | x_2 |  }\right)^{1/2}.
%  $$
%  Then  $0 \le i < b_{j+1}$. Put $s_{0,j}= s_j$ when $i=0$ and $s_{i,j} = s_{j-1} +i s_j$ when  $i\ge 1$. 
% By construction, we have therefore the estimate 
% $$
 % \left( {| y_2 |  q_{k}\over | x_2 |  }\right)^{1/2} \ge 
%  s_{i,j}  \ge 
%  \max  \left( s_j,   \left( {| y_2 |  q_{k}\over | x_2 |  }\right)^{1/2} - s_j \right) \ge {1\over 2} \left( {| y_2 |  q_{k}\over | x_2 |  }\right)^{1/2}  . 
% \leqno{(5.4)}
% $$
Applying  Lemma 4 and using (5.4),  we obtain  the upper bounds
$$
 | \gamma | \le { 2 | y_2 | \over | x_2 |} q_{k-1}q_k + 4 s_{j} q_k \le
  { 2 | y_2 | \over | x_2 |} q_{k-1}q_k + 4 { | y_2 |^{1/2} \over | x_2 |^{1/2}}q_k^{3/2}
  \leqno{(5.5)}
 $$
 and 
 $$
 \eqalign{
 |  \gamma \bx -\by | \le & { 2 | y_2 | \over s_j s_{j+1}} +  { 5 | x_2 | s_{j} \over q_k}
 \le  \left( {2\over s_j} +5 \right) | x_2 |^{1/2} |y_2 |^{1/2} q_k^{-1/2} 
\cr
 \le & 7 | x_2 |^{1/2} |y_2 |^{1/2} q_k^{-1/2}. 
 \cr}
 \leqno{(5.6)}
 $$
  Recall that   $k$ has been chosen satisfying 
$
q_{k-1} \le  q_k^{1/\om}, 
$
where  $\om > 2$.  Consequently, the  first    term $ (2 | y_2 | /  | x_2 | )  q_{k-1}q_k $ occurring on the right hand side of (5.5) is much smaller
than the second one, as  $k$ tends  to infinity. The upper bound 
$$
| \gamma |  \le   5  {| y_2 |^{1/2} \over | x_2 |^{1/2} } q_k^{3/2},   \leqno{(5.7)}
$$
is thus valid for  $k$ large enough.  Combining  (5.6) and  (5.7), we obtain 
$$
| \gamma \bx -\by | \le 7  \root  3 \of 5 | x_2|^{1/3}| y_2|^{2/3} |\gamma|^{-1/3} =  c' | \gamma |^{-1/3}.
$$
 Note that (5.7) turns out to be an  estimate
$$
| \gamma |  \asymp     {| y_2 |^{1/2} \over | x_2 |^{1/2} } q_k^{3/2},
$$
using (4.3). Hence  the norm of $\gamma$ tends to infinity with $k$, and here again, our construction furnishes  infinitely many solutions of 
the inequation (1.4). 

The  assertion (1.4) of Theorem 1 is finally established for any point $\by$ with irrational slope.

\section{5.3. Rational slopes}
We   deduce from Lemma 5  the claim (1.3) of Theorem 1. For any  large integer $k$, it furnishes  a matrix $\gamma \in \Gamma$
satisfying the inequalities
$$
| \gamma | \le  3 { | y_2 | \over | x_2 |} q_{k-1}q_k  \le 3 { | y_2 | \over | x_2 |} q_k^2
\and
| \gamma \bx -\by |  \le { 2 b |x_2 |\over  q_k },
$$
which imply
$$
| \gamma \bx -\by | \le 2 \sqrt{3}b |x_2 |^{1/2} | y_2|^{1/2} | \gamma |^{-1/2} = c  | \gamma |^{-1/2}.
$$
Using  the lower bound for  $\gamma$ given  in Lemma 5, we find the estimate
$$
| \gamma | \asymp { | y_2| \over | x_2 |} q_{k-1}q_k.
$$
Therefore, our construction produces infinitely many solutions $\gamma$ of the inequation (1.3). 
The proof of Theorem 1 is complete.

\section{6. Lower bounds of exponents}

Applying further lemmas 4 and 5,  we now  estimate from below   the exponents $\mu(\bx,\by)$ and $\hat\mu(\bx,\by)$.

\section
{ 6.1. Lower bounds for irrational slopes}
We assume here that the slope $y$ of the point $\by$ is an  irrational number.
As an immediate consequence of (1.4), we get  the lower bound $\mu(\bx,\by)\ge 1/3$.

We prove in this section the lower bound 
$$
\hat\mu(\bx,\by) \ge { \om(y) +1\over 2( 2\om(y)+1)\om(\xi)}, 
$$
claimed  in (1.8). The irrationality measure $\om(y)$ of the slope of the point $\by$ is taken into account
thanks  to the following

\proclaim
Lemma 6. Set 
$$
\tau ={\om(y) \over 2 \om(y) +1}. 
$$
For any $\varepsilon >0$ and any integer $k$ sufficiently large
in terms  of $\varepsilon$, there exists $\gamma\in \Gamma$ such that
$$
| \gamma  | \le C q_k^2 \and | \gamma \bx - \by | \le q_k^{\tau -1 + \varepsilon},
$$
where $C=C(\bx, \by, \varepsilon)$ does not depend upon $k$. 

\pro Once again, we apply  Lemma 4. Let $j$ be the integer defined  by the inequality
$$
s_j \le q_k^\tau < s_{j+1}. \leqno{(6.1)}
$$
Observe that $1/3 \le \tau \le 1/2$, since $\om(y) \ge 1$. Therefore $j$ tends to infinity, as $k$ tends to infinity.
When $\om(y)$ is finite, the lower bound $s_j \ge s_{j+1}^{1/\om}$ holds  for 
any $\om > \om(y)$ provided that $j$ is large enough. Selecting properly 
$\om$ close to $\om(y)$, it follows from (6.1) that
$$
s_j \ge q_k^{\tau /\om(y) -\varepsilon}, \leqno{(6.2)}
$$
for all sufficiently large integers $k$. When $\om(y) = + \infty$, we  read (6.2) as the obvious lower bound
$s_j \ge q_k^{-\varepsilon}$.
Now, Lemma 4 provides  a matrix $\gamma\in \Gamma$ satisfying
$$
| \gamma | \ll q_{k-1}q_k + s_jq_k \le C q_k^2,
$$
and 
$$
| \gamma \bx - \by | \ll {1\over s_j s_{j+1}} + {s_j \over q_k} \ll q_k^{-\tau - \tau/\om(y) + \varepsilon} + q_k^{\tau -1},
$$
by   (6.1) and (6.2). Note that the exponents $ -\tau - \tau/\om(y)$ and $\tau -1$ arising above, are equal by the definition  of $\tau$.
Therefore, we obtain the bound
$$
| \gamma \bx - \by | \ll q_k^{\tau -1 + \varepsilon},
$$
and, decreasing possibly $\varepsilon$, Lemma 6 is proved. \cqfd

\bigskip
For any real number $T$ sufficiently large, let $k$ be the integer defined  by the inequalities
$$
C q_k^2 \le T < Cq_{k+1}^2. \leqno{(6.3)}
$$
Clearly, $k$ tends to infinity when $T$ tends to infinity. For any $\varepsilon >0$, we can bound further
$$
T \le C q_{k+1}^2 \le q_k^{2 \om(\xi)+ \varepsilon}, \leqno{(6.4)}
$$
when $T$ is large enough. Then,  Lemma 6  gives a  matrix  $\gamma \in \Gamma$ satisfying
$$
| \gamma  |  \le C q_k^2 \le T \and | \gamma \bx - \by  | \le q_k^{\tau -1 + \varepsilon} \le T^{-(1- \tau - \varepsilon)/(2 \om(\xi) + \varepsilon)},
$$ 
by (6.3) and (6.4). Therefore
$$
\hat\mu(\bx,\by) \ge { 1 - \tau - \varepsilon \over 2 \om(\xi) + \varepsilon},
$$
and letting $\varepsilon$ tends to $0$, we obtain the claimed   lower bound 
$$
\hat\mu(\bx,\by) \ge {1- \tau \over 2 \om(\xi)} =  { \om(y) +1\over 2( 2\om(y)+1)\om(\xi)}. 
$$

\section
{   6.2. Lower bounds for rational slopes}
In this section, we prove   that the lower  bounds
$$
\mu(\bx,\by) \ge { \om(\xi)\over \om(\xi) + 1} \and \hat\mu(\bx,\by) \ge { 1 \over \om(\xi) +1}
$$
hold for any  point $\by $ with rational slope $y$, or when $y_2=0$.
As  in Section 4.2, we assume that $y_2\not= 0$ and that  
   $y = a/b$,  where $a$ and $b$ are coprime integers 
with $| a | \le b$. 

We start with the inequality $\mu(\bx,\by) \ge \om(\xi)/(\om(\xi) +1)$. For any $\om < \om(\xi)$ there exist infinitely many integers $k$
satisfying  $q_{k-1}\le q_k^{1/\om}$. Using Lemma 5 for such an index $k$, we get $\gamma\in \Gamma$ such that
$$
 | \gamma | \ll q_{k-1}q_k \ll q_k^{1 + 1/\om} \and | \gamma \bx - \by | \ll q_k^{-1}.
 $$
 Then  $| \gamma \bx - \by | \ll | \gamma |^{-\om/(\om +1)}$ for infinitely many $\gamma$. Hence $\mu(\bx, \by)\ge \om(\xi)/(\om(\xi) + 1)$
 by letting $\om $ tend to $\om(\xi)$.

As for the lower bound $\hat\mu(\bx,\by) \ge 1/(\om(\xi)+1)$, we briefly take again  the  argumentation given  in Section 6.1.
We may obviously assume that $\om(\xi)$ is finite.  
For any real number $T$ sufficiently large, let $k$ be the integer uniqueley determined by 
$$
 3 { | y_2 | \over | x_2 |}  q_{k-1} q_k \le T <   3 { | y_2 | \over | x_2 |}  q_k q_{k+1}.
$$
For any $\varepsilon >0$, we bound from above
$$
T \le  3 { | y_2 | \over | x_2 |} q_k  q_{k+1} \le q_k^{ \om(\xi)+ 1+  \varepsilon}, 
$$
when $k$ is large enough. Lemma 5 gives us  a matrix  $\gamma \in \Gamma$ satisfying
$$
| \gamma  |  \le 3 { | y_2 | \over | x_2 |}  q_{k-1} q_k  \le T 
\and
 | \gamma \bx - \by  | \le  { 2 b |x_2 |\over  q_k } \le  2 b |x_2 | T^{-1/( \om(\xi) + 1 + \varepsilon)}.
$$ 
 Therefore  $ \hat\mu(\bx,\by) \ge 1/ ( \om(\xi) +1 +  \varepsilon)$ for any $\varepsilon>0$.

\section
{7. Proof of Theorem 3}

Recall the matrices $M_{k}$ and $N_j$ introduced in  Sections  2 and 4.1.
We intend to show that if an element $\gamma\bx$ of the orbit  is close to a given point $\by$, then $\gamma$ can be factorized
in the form  $\gamma = N_j G M_{k}$, with a good  estimate of  the norm $| G|$ for suitable indices $j$ and $k$. 
Without loss of generality, we may assume that   $\bx = \left(\matrix{ \xi  \cr 1 \cr}\right)$. 

\proclaim
Lemma 7. Let  $k$ be a positive integer, $\mu$ and $T$ be   real numbers  such that    
$$
0  \le  \mu \le 1 \and  q_{k-1}q_k \le T \le q_kq_{k+1},
$$  
 and let  $\gamma \in \Gamma$ satisfy
$$
 |  \gamma | \le 2 T   \and | \gamma \bx - \by | \le T^{-\mu  } .
$$
Let  $j$ be a positive integer such that 
$
s_j \ge T^{\mu/2}.
$
Then   $ \gamma $ can be decomposed as a product  
$
\gamma = N_j G M_{k}, 
$
 where  the two columns of the matrix  $G = \left(\matrix{ m & \ell \cr m'  & \ell'\cr}\right) \in Ê\Gamma $ 
 satisfy the norm bound 
$$
\max ( | m |, |m' | )  \le { c  s_{j}  T^{1 -\mu}\over q_k}
\and 
\max ( | \ell  |, |\ell'  | )  \le c s_{j} q_k  T^{-\mu},
$$
with  $c = 10 \max( | \by | , | \by |^{-1})$.

\pro
Write 
$
\gamma =  \left(\matrix{ v_1 & u_1\cr v_2& u_2\cr}\right)
$
 and put 
$$
\Lambda_1 = v_1 \xi +  u_1 -y_1 , \quad  \Lambda_2 = v_2 \xi +  u_2 -y_2.
$$
 The upper bound  $  | \gamma \bx - \by | \le   T^{-\mu} $  means that  
$$
\max( | \Lambda_1|, |\Lambda_2| ) \le T^{-\mu} . \leqno{(7.1)}
$$
 We have  the identities
$$
\eqalign{
v_1 y_2 - v_2y_1 = & \left|\matrix{ v_1 & y_1\cr v_2  & y_2\cr}\right|= 
 \left|\matrix{ v_1 & v_1 \xi +u_1-\Lambda_1 \cr v_2 & v_2 \xi + u_2-\Lambda_2 \cr}\right| 
=  1 + \Lambda_1 v_2 - \Lambda_2 v_1,
\cr 
u_1 y_2 - u_2y_1 = & \left|\matrix{ u_1 & y_1\cr u_2  & y_2\cr}\right|= 
 \left|\matrix{ u_1 & v_1 \xi +u_1-\Lambda_1 \cr u_2 & v_2 \xi + u_2-\Lambda_2 \cr}\right| 
=  -\xi + \Lambda_1 u_2 - \Lambda_2 u_1. 
\cr }
\leqno{(7.2)}
$$
By (7.1), they imply the upper bound
$$
\max \Big( | u_1 y_2 - u_2 y_1 |,  | v_1 y_2 - v_2 y_1 |\Big)  \le 1 + 4 T^{1-\mu}. \leqno{(7.3)}
$$

We first factorize  $N_j$. Define 
$$
\eqalign{
\gamma ' = & N_j^{-1}\gamma = \left(\matrix{ t_j & t_j'\cr s_j& s_j'\cr}\right)^{-1} \left(\matrix{ v_1 & u_1\cr v_2& u_2\cr}\right)
=  \left(\matrix{ s_j'v_1 -t_j'v_2 & s_j'u_1-t_j'u_2\cr -s_jv_1+ t_jv_2& -s_ju_1+ t_ju_2\cr}\right)
\cr 
 = & 
 {1\over y_2}   \left(\matrix{   s_j'(v_1y_2 - v_2y_1) +v_2(s_j'y_1-t_j'y_2) &  s_j'(u_1y_2 - u_2y_1) +u_2(s_j'y_1-t_j'y_2)
   \cr
   - s_j(v_1y_2 - v_2y_1) -v_2(s_jy_1-t_jy_2)   &  -s_j(u_1y_2 - u_2y_1) -u_2(s_jy_1-t_jy_2) \cr}\right).
 \cr}
$$
Using  (7.3) and the estimate $ |s_jy-t_j| \le |s_j'y-t_j'|\le 1/s_j$, we deduce from the above expression   the upper bound for the norm
$$
| \gamma' |  \le { s_j(1+ 4T^{1-\mu})\over | y_2|}+ {2 T \over s_j} \le (5 |y_2 |^{-1}+2)  s_j T^{1-\mu},\leqno{(7.4)}
$$
since $s_j\ge T^{\mu/2}$. 
Now, put $\gamma' =  \left(\matrix{ v'_1 & u'_1\cr v'_2& u'_2\cr}\right)$ and write 
$$
 \left(\matrix{v'_1\xi + u'_1 \cr v'_2 \xi +u'_2 \cr}\right) = \gamma'\bx = N_j^{-1}\gamma \bx = N_j^{-1} \left(\matrix{ y_1 +\Lambda_1\cr y_2 +\Lambda_2\cr}\right)
 = \left(\matrix{y_1s'_j -y_2t'_j + s'_j\Lambda_1 - t'_j\Lambda_2 \cr- y_1s_j +y_2t_j -s_j\Lambda_1 + t_j\Lambda_2 \cr}\right).
$$
It follows that
$$
\max \Big( | v'_1 \xi+ u'_1|, | v'_2 \xi +u'_2 | \Big) = | \gamma'\bx |  \le { | y_2 | \over s_j} + 2 s_j T^{-\mu}
 \le ( |y_2| +2) s_j T^{-\mu}.\leqno{(7.5)}
$$
Now, we multiply $\gamma'$ on the right by $M_{k}^{-1}$ and set
$$
G  = N_j^{-1} \gamma M_{k}^{-1} = \gamma'M_{k}^{-1}.
$$
Suppose first that $k$ is even. We find the formula
$$
G = \left(\matrix{ v'_1 & u'_1\cr v'_2& u'_2\cr}\right) \left(\matrix{ q_{k} & -p_{k}\cr -q_{k-1}& p_{k-1}\cr}\right)^{-1}
 =   \left(\matrix{  p_{k-1}v'_1 + q_{k-1} u'_1&  p_{k}v'_1  + q_{k} u'_1\cr   p_{k-1}v'_2 + q_{k-1} u'_2&   p_{k}v'_2 + q_{k} u'_2\cr}\right).
 $$
Write next
$$
\eqalign{
 \ell = & p_{k}v'_1 + q_{k} u'_1    = - v'_1( q_{k}\xi - p_{k})   + q_{k}( v'_1\xi +u'_1) ,
 \cr
 \ell'= & p_{k}v'_2 + q_{k} u'_2    = - v'_2( q_{k}\xi - p_{k})   + q_{k}( v'_2\xi +u'_2) ,
 \cr
m= & p_{k-1}v'_1 + q_{k-1} u'_1    = - v'_1( q_{k-1}\xi - p_{k-1})   + q_{k-1}( v'_1\xi +u'_1) ,
 \cr
m'= &  p_{k-1}v'_2 + q_{k-1} u_2    = - v'_2( q_{k-1}\xi - p_{k-1})   + q_{k-1}( v'_2\xi +u'_2) .
 \cr}
$$
We deduce from (2.1),  (7.4) and  (7.5) that 
$$
\displaylines{
\max ( | \ell |, |  \ell' |  )  \le (5 |y_2 |^{-1}+2)  {s_j T^{1-\mu}\over q_{k+1}}  + ( |y_2| +2) q_{k}s_j T^{-\mu} 
\le c  s_{j}q_k T^{-\mu},
\cr
\max ( | m |, | m' |  )  \le (5 |y_2 |^{-1}+2)  {s_j T^{1-\mu}\over q_{k}}  + ( |y_2| +2) q_{k-1}s_j T ^{-\mu} 
\le {c s_j T^{1-\mu}\over q_{k}} ,
\cr}
$$
since $q_{k-1}q_k \le T \le q_kq_{k+1}$.  The case $k$ odd  leads to the same upper bound. \cqfd

 \bigskip
 
We are now able to prove  Theorem 3. 
 Let  $\cC$ be a   compact subset  of the punctered line   $(\bR\setminus \{0\}) \left(\matrix{ y  \cr 1 \cr}\right)$, 
 and let $\mu$ be a real number greater than $1/2$. 
 Denote by $\cC_\mu$ the subset 
 consisting of    the points  $\by\in \cC$ for which the inequation 
 $$
  | \gamma \bx - \by | \le    | \gamma|^{-\mu} \leqno{(7.6)}
  $$
  has infinitely many solutions  $\gamma \in \Gamma$.
  We have to  show that $\cC_\mu$ has null Lebesgue measure. 
  
  Let  $\gamma \in\Gamma $ and  $\by \in \cC_\mu$ satisfying (7.6). 
  Assuming that $|\gamma |$ is large enough, let $k\ge 1 $ and  $n\ge 0 $ be the integers    defined  by the inequalities
  $$
  q_{k-1}q_k  < | \gamma | \le q_kq_{k+1} \and 2^n q_{k-1}q_k <  | \gamma |  \le 2^{n+1} q_{k-1}q_k. \leqno{(7.7)}
  $$
 Put $T = 2^nq_{k-1}q_k$. It follows from   (7.6) and (7.7) that 
    $$
  | \gamma | \le 2 T  \and  | \gamma \bx - \by | \le  | \gamma | ^{-\mu} \le T^{-\mu}. \leqno{(7.8)}
  $$
  Let $j$ be the smallest integer such that $s_j \ge T^{\mu/2}$. Since we have assumed that $\om(y)=1$, for any positive
  real number $\varepsilon$, we can bound from above $s_j \le T^{\mu/2  + \varepsilon}$ when $j$ is large enough. 
  Note that $j$ is arbitrarily large if we take $\gamma$ of sufficiently large norm. 
 Then,  Lemma 7 provides us with the decomposition 
  $
  \gamma = N_j G M_{k}
    $
  for some matrix  $G = \left(\matrix{ m & \ell \cr m'  & \ell'\cr}\right)$ in  $\Gamma$ whose columns satisfy the bound of norm  
  $$
  \eqalign{
  \max ( | m |, |m' | )  \le &  {c s_j T^{1-\mu}\over q_{k}}  \le {c T^{1-\mu/2 +\varepsilon}\over q_k} =  B_1,
\cr
\max ( | \ell  |, |\ell'  | )  \le  &  c  s_{j}q_k T^{-\mu} \le c q_k T^{-\mu/2 + \varepsilon} =  B_2,
\cr} \leqno{(7.9)}
$$
where  the coefficient $c=10\max_{\by\in \cC}(| \by |, |\by |^{-1})$ depends only upon  $\cC$. 

It is easily seen  that the set of matrices $G\in \Gamma$ whose  first and second columns  have   norm respectively bounded by $B_1$ and $B_2$, has
at most $ 4(2B_1+1)(2B_2+1)$ elements. Of course, if either  $B_1$ or $B_2$ is smaller than $1$, no such matrix exists. 
Hence,  there are  at most 
$$
36  B_1 B_2 = 36  c^2 T^{1-\mu + 2 \varepsilon}
$$
 matrices $G$ in $\Gamma$ satisfying (7.9).
The second upper bound  of (7.8) means   that  $\by$ belongs to the intersection  of  the line $\bR \left(\matrix{ y  \cr 1 \cr}\right)$
with  the square centered at  the point $N_jGM_k\bx$ with  side $ 2  T^{-\mu}  $. This intersection   is  a 
  segment  of  Euclidean length  $\le 2  \sqrt{2}T^{-\mu}$. 
 For fixed  $k$ and $n$,  at most  $36 B_1 B_2 $ such  segments   may thus appear. 
 It follows that  $\by$ belongs to some subset of the  line  $\bR \left(\matrix{ y  \cr 1 \cr}\right)$ whose  Lebesgue measure does not exceed
 $$
(36 B_1 B_2)( 2 \sqrt{2} T^{-\mu}) =  72 \sqrt{2} c^2 (2^nq_{k-1}q_k)^{1-2\mu+ 2\varepsilon}.
 $$
 
Note that  the sequence  $q_k$ of denominators of  convergents of the irrational number $\xi$ is bounded from below by the 
Fibonacci sequence $1,1,2,\ldots$, which  grows geometrically. 
Therefore,   the series
$$
\sum_{k\ge 1}  \sum_{ n\ge 0}   (2^nq_{k-1}q_k)^{1-2\mu+ 2\varepsilon}
$$
converges when $\varepsilon$ is small enough, since the exponent $1  -2 \mu$ is negative.
 By Borel-Cantelli Lemma, the $\limsup$ set $\cC_\mu$  has  null Lebesgue measure.

 \section{8. Upper bounds for rational slopes}
 Here we  prove  that the upper bounds  
 $$
\mu(\bx,\by) \le { \om(\xi)\over \om(\xi) + 1} \and \hat\mu(\bx,\by) \le { 1 \over \om(\xi) +1}
$$
hold  for  any  point  $\by\not= {\bf 0}$  with  rational slope $y$.
 Since  the reverse inequalities have been  established in Section 6.2,  the proof of (1.7) will then be complete.
  To that purpose, we adapt to  rational slopes  the factorisation
 method displayed in the preceding section. We obtain  the following explicit lower bound of distance  which may have its own  interest.
 
 \proclaim
 Theorem 4.  Let $\by =  \left(\matrix{y_1\cr y_2\cr}\right)$ be a point having 
  rational slope $y_1/y_2 = a/b$, where  $a$ and $b$ are coprime integers with $| a | \le b$,  and 
 let $k$ be a positive integer such that  $ q_k \ge 12 b/| y_2 | $. 
 Then, for any $\gamma \in \Gamma$ with norm 
 $$
  \left\vert \gamma\right\vert \le { | y_2 | \over 4} q_k q_{k+1},
  $$
     we have the lower bound
 $$
 \left\vert  \gamma \left(\matrix{\xi\cr 1\cr}\right) - \by \right\vert \ge { 1\over 4 b  q_k}.
 $$
 
\pro
Recall the  matrix $N=  \left(\matrix{ a &a'\cr b& b'\cr}\right)$ in $ \Gamma$
introduced in Section 4.2. 
 Notice that
$N^{-1}$  maps the line $ \bR  \left(\matrix{a\cr b\cr}\right)$ on to the horizontal axis $\bR\left(\matrix{1\cr 0\cr}\right)$. 
Therefore any point close to the line $ \bR  \left(\matrix{a\cr b\cr}\right)$ is sent by the map $N^{-1}$ to a point close to the horizontal axis.
We insert this additional information into  the proof of Lemma 7 with   $\mu =1/2$. 

 Set 
$
 \left(\matrix{\Lambda_1 \cr \Lambda_2 \cr}\right)=  \gamma \left(\matrix{\xi\cr 1\cr}\right) - \by 
$
and suppose on the contrary that 
 $
\max(|\Lambda_1| , |\Lambda_2|)<  (4 b q_k)^{-1}.
 $
 Put
$$
\gamma ' = N^{-1}\gamma =  \left(\matrix{ v'_1 & u'_1\cr v'_2& u'_2\cr}\right). 
$$
Noting that 
$$
b y_1 - a y_2 = 0  \and b'  y_1 - a' y_2 = {y_2\over b},
$$
we obtain as in Section 7 the expressions 
$$
\gamma' =  \left(\matrix { \displaystyle
 {b'(v_1y_2 -v_2y_1)\over y_2} +{v_2\over b} 
&  \displaystyle
 { b'(u_1y_2 -u_2y_1)\over y_2} +{u_2\over b}
\cr  \displaystyle
 -{ b(v_1 y_2 -v_2y_1)\over y_2} 
&  \displaystyle
   -{ b(u_1 y_2 -u_2y_1)\over y_2}
\cr}\right)
\leqno{(8.1)}
$$
and 
$$
\gamma'\bx = \left(\matrix{  {y_2\over b}  + b'\Lambda_1 - a'\Lambda_2 \cr -b\Lambda_1 + a\Lambda_2 \cr}\right). 
\leqno{(8.2)}
$$
Using the formulas (7.2), we have that 
$$
  | v_1 y_2 - v_2 y_1 |   \le 1 +2\max( | \Lambda_1 |, | \Lambda_2 |)  | \gamma  | 
  \le  1+ { | y_2 | \over 8 b} q_{k+1}  \le { | y_2 | \over 4 b} q_{k+1} , 
\leqno{(8.3)}
$$
since we  have assumed that $q_k \ge 12 b/|y_2 |  $. Then,  we deduce   from   the expressions (8.1), (8.2) and from the upper bound (8.3)  that 
$$
  |v'_2| <  {q_{k+1} \over 4}  \and  | v'_2\xi + u'_2 | <  {1  \over 2 q_k}.  \leqno{(8.4)}
$$
Set now 
$$
G  = N^{-1} \gamma M_k^{-1} = \gamma'M_k^{-1}.
$$
Assuming that $k$ is even (the case $k$ odd is similar), we use again the expressions
$$
G =
   \left(\matrix{ - v'_1( q_{k-1}\xi - p_{k-1})   + q_{k-1}( v'_1\xi +u'_1) & - v'_1( q_{k}\xi - p_{k})   + q_{k}( v'_1\xi +u'_1) \cr   
  - v'_2( q_{k-1}\xi - p_{k-1})   + q_{k-1}( v'_2\xi +u'_2) & - v'_2( q_{k}\xi - p_{k})   + q_{k}( v'_2\xi +u'_2)  \cr}\right)
 $$
 obtained in Section 7. We deduce from (2.1) and (8.4) the upper bound
 $$
  \Big\vert- v'_2( q_{k}\xi - p_{k})   + q_{k}( v'_2\xi +u'_2) \Big\vert
  \le {   |v'_2|\over q_{k+1} } + q_k   | v'_2\xi + u'_2 | \le  { 1\over 4 } + {1\over 2} <1 ,
$$
for the absolute value of the lower right entry  of the matrix $G$,  which therefore  vanishes.
It follows that $G$ has the form
$$
G = \pm  \left(\matrix{ m  & -1 \cr   
 1  & 0 \cr}\right),
 $$
where $m$ is an integer. Hence
$$ 
 \left(\matrix{{y_2\over b}  + b'\Lambda_1 - a'\Lambda_2 \cr -b\Lambda_1 + a\Lambda_2 \cr}\right)
 = \gamma' \bx = G M_k\bx  = \pm  \left(\matrix{ m \epsilon_k - | \epsilon_{k-1}|   \cr   \epsilon_k \cr}\right).
$$
Looking at the first component of the above vectorial equality, we find the estimates
$$
{ | y_2 | \over b} -{1\over 2 q_k} \le 
\Big\vert {y_2\over b}  + b'\Lambda_1 - a'\Lambda_2  \Big\vert 
=
\Big\vert   m \epsilon_k - | \epsilon_{k-1}| \Big\vert  \le { | m | \over q_{k+1} } +{ 1\over q_k}. 
$$
We thus obtain the lower bound 
 $$
 | m  |  \ge { |y_2 | q_{k+1}\over 2 b} \ge 6,\leqno{(8.5)}
 $$
 since $q_k \ge 12 b/ | y_2| $. Now, write 
 $$
 \eqalign{
 \gamma = &\pm  \left(\matrix{ a   & a' \cr   b  & b' \cr}\right) \left(\matrix{ m  & -1 \cr    1  & 0 \cr}\right)
  \left(\matrix{ q_k  & -p_k \cr  - q_{k-1}  &  p_{k-1} \cr}\right)
  \cr
  = & \pm  \left(\matrix{ a m q_k  + a q_{k-1} +a'q_k   & -amp_{k-1} - ap_{k-1} -a'p_{k}
   \cr  b m q_k +  b q_{k-1} +b'q_k   & -bmp_{k-1}  - bp_{k-1} -b'p_{k} \cr}\right).
   \cr}
   $$
 Hence,   taking (8.5) into account, we find the lower bound 
 $$
 | \gamma | \ge b( | m | -2) q_k \ge { | y_2 | \over 3  } q_k q_{k+1},
 $$
 which  contradicts  the  assumption $ | \gamma | \le ( | y_2 | / 4  ) q_k q_{k+1}$.
 \cqfd

\medskip 
 We first deduce from Theorem 4 that 
 $
 \mu(\bx, \by ) \le \om(\xi)/(\om(\xi)+1) .
 $
For any matrix  $\gamma $  in $\Gamma$ with norm  $|\gamma |$ large enough, 
 let $k$ be the integer defined by the inequality 
  $$
    {| y_2 | \over 4 }q_{k-1}q_{k}< | \gamma | \le {| y_2 | \over 4 }q_kq_{k+1}.
  $$
  In the case where $\om(\xi)$ is finite, let  $\om $ be a real number greater  than  $\om(\xi)$. 
  We then bound from below    $q_{k-1} \ge q_k^{1/\om}$,  
if $k$ is large enough in terms  of  $\om $.   In the case  $ \om(\xi)= +\infty$, we simply bound from below $q_{k-1} \ge 1$. 
Now, Theorem 4 gives us the lower bound
$$
 | \gamma \bx - \by | \ge   { 1 \over 4b q_k} \ge {1 \over 4b} {1 \over  (4  | \gamma |/| y_2 |)^{1/(1+1/\om)}} ,
 $$
 where the exponent $1/(1+1/\om)$ is understood to be 1 when $\om(\xi)= +\infty$. 
  The latter  lower bound of distance is valid for any $\gamma \in \Gamma$ with large norm. It thus implies   the upper bound 
 $$
  \mu(\bx,\by)  \le {1 \over 1+ {1 \over \om}}= {\om\over \om+1}.
  $$
  Letting  $\om $ tend to $\om(\xi)$, we have proved  the claim.
  
 Let $\mu$ be  a positive real number such that the inequations 
  $$
  |  \gamma | \le T \and | \gamma\bx - \by | \le T^{-\mu} \leqno{(8.6)}
  $$
  have  a solution $\gamma \in \Gamma$ for any  large real number $T$. 
  Let $\om$ be a real number smaller than $\om(\xi)$. There exist infinitely many integer $k$ such that $q_{k+1} \ge q_k^\om$.
   Choose   $T=  (| y_2 |/ 4 )q_kq_{k+1}$ for  such an integer $k$. Thus 
     $T \ge (| y_2 |/ 4 )q_k^{1+\om}$, and Theorem 4 now gives  the lower bound
$$
 | \gamma \bx - \by | \ge   { 1 \over 4b q_k} \ge {1 \over 4b} {1 \over  (4 T /| y_2 |)^{1/(1+\om)}} ,
 $$
for any $\gamma \in \Gamma $ with norm $| \gamma | \le T$. Comparing with (8.6), we find that $\mu \le 1/(1+\om)$.  
Letting $\om$ tend to $\om(\xi)$, we obtain the expected bound  
$\hat\mu(\bx,\by) \le 1/(\om(\xi) +1)$. 

 \section
{9. Approximation with signs}

Let us first state   a  theorem due to  Davenport  and   Heilbronn which gives a version of  Minkowski Theorem with prescribed signs
\cite{DaHe}. 

\proclaim 
Theorem {\rm (Davenport--Heilbronn)}. Let  $\xi$ be an  irrational number  and let  $y$ be a real number not belonging to the subgroup  $\bZ\xi + \bZ$.
 There exist  infinitely many pairs of integers $(v,u)$ such that  
 $$
v > 0 \and  0 < v\xi + u -y  \le  {1\over v}.
 $$
  
Here is an analogous statement  for  $\Gamma$-orbits. For simplicity, we assume that the  target point  $\by= \left(\matrix{y_1 \cr y_2  \cr}\right)$
belongs to the positive  quadrant  $\bR_+^2$. 

\proclaim
Theorem 5. Let   $\xi $ be an irrational number and let  $y_1$, $y_2$ be two positive real numbers
such  that  the ratio  $y=y_1/y_2$  is an irrational number with irrationality measure
 $\om(y) = 1$.  Then,    for any positive real number  $\mu < 1/3$, there exist infinitely many 
matrices $\gamma =  \left(\matrix{ v_1 & u_1\cr v_2 & u_2 \cr}\right) \in \Gamma$  satisfying 
$$
v_1 > 0, \quad v_2  >0  \and   0 < v_1   \xi + u_1 -y_1  \le | \gamma |^{-\mu},
 \quad  0 < v_2   \xi + u_2 -y_2  \le | \gamma |^{-\mu} . 
   $$

\noindent
{\bf Remark.}  Other constraints  of signs are admissible.
 Notice however that  $v_1$ and  $v_2$ have necessarily the same sign whenever  $y_1$ and $y_2$ have the  same sign, 
if  we assume that  
$\left\vert \gamma   \left(\matrix{\xi \cr 1  \cr}\right) - \left(\matrix{y_1 \cr y_2  \cr}\right)\right\vert = \cO( | \gamma|^{-\mu})$ with $\mu > 0$. 
That follows from the estimate  
$$
\eqalign{
v_1y_2 -v_2y_1 =  \left|\matrix{ v_1 & y_1\cr v_2  & y_2\cr}\right|= & \left|\matrix{ v_1 & v_1 \xi +u_1 \cr v_2 & v_2 \xi + u_2 \cr}\right| 
-  \left|\matrix{ v_1 & v_1   \xi + u_1 -y_1\cr v_2  &v_2   \xi + u_2 -y_2 \cr}\right|\cr
= & Ê\, Ê\, 1 +  \cO\left( | \gamma |^{1-\mu}\right),
\cr}
$$
already mentioned  in (7.2).  Theorem 5 is   a sample of  statements that could be obtained by reworking the 
previous sections and   controling  all   signs. 

\medskip

Denote by $\Gamma_+$ the semi-group of $\Gamma$ consisting of the matrices $\gamma$ with non-negative entries. 
Theorem 5 enables us to recover in a constructive way the following result from  \cite{DaNo}: 
\proclaim
Corollary {\rm (Dani-Nogueira)}. Let $\xi$ be a negative irrational number. Then,   
the intersection with $\bR_+^2$ of the semi-orbit $\Gamma_+ \left(\matrix{ \xi \cr 1\cr}\right)$ is  dense in $\bR_+^2$.

\pro  The  points $\by = \left(\matrix{y_1 \cr y_2  \cr}\right)\in \bR_+^2$ for which the slope $y=y_1/y_2$ 
has irrationality measure $\om(y)=1$  form a full set in  $\bR_+^2$
({\it i.e.} the complementary set has null Lebesgue measure), hence dense. For any such point  $\by$,  Theorem 5 
provides us with a sequence of points in  $\Gamma_+   \left(\matrix{\xi \cr 1  \cr}\right)$ tending to  $\by$, since  the second column
 $ \left(\matrix{  u_1\cr  u_2 \cr}\right)$  of $\gamma$ has necessarily positive entries when 
  $v_1 > 0$, $v_2 >0$, $\xi <0$    and  $ \left\vert \gamma  \left(\matrix{\xi \cr 1  \cr}\right) - \by \right\vert $  is sufficiently small. \cqfd

\bigskip
\noindent
{  \bf Proof of Theorem 5.}
We take again the construction  of Section  4.1.  In order to prescribe  positive  signs,   
we need to introduce a variant $ \tilde N_j$ of the  matrices $N_j$ which induces  slight modifications  in  the estimates. 

 Recall that    $(t_j/s_j)_{j\ge 0}$  stands for the sequence of convergents of   $y$.  For any  $j\ge 1$, we set  
 $$
 \tilde N_j =  \left(\matrix{ t_{j-1} &  t_j\cr s_{j-1} & s_j  \cr}\right) \qquad {\rm or  } \qquad 
  \tilde N_j = N_j = \left( \matrix{  t_j &  t_{j-1}\cr  s_j & s_{j-1} \cr}\right) ,   
$$
 respectively when    $j$ is  even  or odd.  
 The matrix  $\tilde N_j$ belongs to  $\Gamma_+$ and has norm  
 $$
 | \tilde N_j |   = \max(s_j,t_j) \asymp s_j.
 $$
Notice that if we put
$$
 \tilde N_j =  \left(\matrix{ t &  t' \cr s & s'  \cr}\right)
 \and 
 \delta = s y -t, \quad \delta' = s'y - t'
 $$
 then  $\delta$ is negative,  and we now have the (weaker) estimates
 $$
 {1\over 2 s_{j+1}} < - \delta  \le {1\over s_j} \and    | \delta' | \le {1\over s_j} \leqno{(9.1)}
 $$
  for any  $j\ge 1$.  We  consider   matrices of the form 
$
\gamma =  \tilde N_j U^{\ell}  M_k  ,
$
where  $k$ and $\ell$ are  positive integers and $k$ is odd.  Observe that the matrix $M_k$ has positive entries on  its first column
precisely when $k$ is odd. We find  the formula
 $$
 \gamma  =   \left(\matrix{ v_1 &  u_1 \cr v_2 & u_2  \cr}\right) = \left(\matrix{  \ell t q_{k-1} +  t q_k + t'q_{k-1} &  -\ell t p_{k-1} -  t p_k - t'p_{k-1} 
\cr   \ell s q_{k-1} + s q_k + s'q_{k-1}  & -\ell s p_{k-1} - s p_k - s'p_{k-1} \cr}\right) .
$$
It follows that the first column  $ \left(\matrix{v_1 \cr v_2\cr}\right)$ of the matrix  $\gamma$ has positive entries, and that 
we have the bound of  norm  
$$
| \gamma |  \le ( \ell + 2) | \tilde N_j  |  |  M_k |  \ll \ell s_j q_k.   \leqno{(9.2)}
$$
 Denote as usual 
 $$
  \left(\matrix{ \Lambda_1\cr \Lambda_2 \cr}\right) =  \left(\matrix{ v_1\xi +u_1 -y_1\cr v_2 \xi +u_2 -y_2\cr}\right) .
 $$
 Taking again the computations of Lemma 3, we  find the formulas
 $$
 \Lambda_1 - y \Lambda_2 = -\delta(\epsilon_k + \ell \vert\epsilon_{k-1}\vert) - \delta'\vert\epsilon_{k-1}\vert \leqno{(9.3)}
 $$
 and 
 $$
 \Lambda_2 = s\vert\epsilon_{k-1}\vert(\ell - \rho) \quad {\rm with} \quad \rho = 
 { y_2 \over s \vert\epsilon_{k-1}\vert} - { \epsilon_k \over \vert\epsilon_{k-1}\vert} - {s' \over s}.
 \leqno{(9.4)}
 $$

For any odd large index $k$, let $j$ be the integer defined  by the estimate
$$
s_{j-1} <  q_k^{1/3} \le s_j.
$$
Since we have assumed that $\om(y)=1$, the inequalities
$$
q_k^{1/3 - \varepsilon} \le s_{j-1}  < q_k^{1/3}  \le s_j \le q_k^{ 1/3 + \varepsilon} \and  s_{j+1} \le q_k^{1/3 + 2 \varepsilon}\leqno{(9.5)}
$$
hold for any $ \varepsilon >0$, provided that $j$ is large enough.  
We deduce from the expression for $\rho$,  given in (9.4),  the estimate
$$
 y_2 q_k^{2/3-\varepsilon} - 1 - q_k^{2\varepsilon} \le \rho \le 2 y_2 q_k^{2/3 + \varepsilon } +1 , \leqno{(9.6)}
$$
using  (2.1), (9.5),   and noting that $ 0 \le s'/s \le s_j/s_{j-1} \le  q_k^{2 \varepsilon}$ by (9.5). It follows that the real
number $\rho$ is positive, when $k$ is large enough. 
Let $\ell$ be the smallest integer larger or equal to  $\rho$. We deduce from  (2.1) and  (9.5) that
$$
0 < \Lambda_2 \le s\vert\epsilon_{k-1}\vert \le  {s_j\over q_k} \le  q_k^{-2/3+ \varepsilon}. \leqno{(9.7)}
$$
  Moreover, $\ell$ is a positive integer satisfying
$$
q_k^{2/3 - \varepsilon}  \ll \ell \ll q_k^{2/3 + \varepsilon}, \leqno{(9.8)}
$$
according to the estimate (9.6). Using (9.5) and (9.8), observe now that the leading term on the right hand side of formula (9.3)
 giving $\Lambda_1-y\Lambda_2$  is $-\delta \ell \vert\epsilon_{k-1}\vert$,  which is positive.  We thus find the estimate
$$
0 < \Lambda_1 - y \Lambda_2 \ll { \ell \vert\epsilon_{k-1}\vert\over s_{j}} \ll q_k^{-2/3 +  \varepsilon}, \leqno{(9.9)}
$$
making use of the inequalities (2.1), (9.1), (9.5) and (9.8). 
Since $y$ is positive, it follows  that $\Lambda_1$ is  positive as well. Moreover, we deduce 
from (9.7) and (9.9) that 
$$
 \max (\Lambda_1, \Lambda_2 )  \ll q_k^{-2/3 + \varepsilon}. \leqno{(9.10)}
$$
Next,   the bound of  norm 
 $$
 | \gamma | \ll \ell s_j q_k \ll q_k^{2 +  2  \varepsilon}.
 $$
follows from (9.5) and (9.8). Now, we deduce  from (9.10) that 
 $$
  \max (\Lambda_1, \Lambda_2 )  \ll | \gamma |^{-(2/3 -  \varepsilon)/( 2 + 2  \varepsilon)} \le | \gamma |^{-\mu},
 $$
 provided  
 $
 \mu  < (2 - 3 \varepsilon)/(  6 + 6 \varepsilon).
 $
 Since  $\mu < 1/3$, this last inequality is   satisfied by choosing  $\varepsilon$ small enough. 
 
 Finally, observe that we have the estimate of norm 
 $$
 | \gamma |  \asymp \ell s q_{k-1} \gg q_k^{1 -2 \varepsilon}q_{k-1},
 $$
 by (9.5) and (9.8). Therefore, $| \gamma | $ may be arbitrarily large when  $k$ is large enough,
  and our construction produces infinitely many matrices  $\gamma$ verifying  Theorem 5. 
  \cqfd

 \vskip 8mm

\centerline{\bf References}

\vskip 5mm

\beginthebibliography{999}

\bibitem{Bu}
Y. Bugeaud,
{\it A note on inhomogeneous Diophantine approximation},
Glasgow Math. J. {\bf 45} (2003), 105--110.

\bibitem{BuLa}
Y. Bugeaud and M. Laurent,
{\it Exponents of inhomogeneous Diophantine approximation},
Moscow Math. J. {\bf 5} (2005), 747--766.

\bibitem{Cas}
J. W. S. Cassels,
An introduction to Diophantine Approximation.
Cambridge Tracts in Math. and Math. Phys., vol. 99, Cambridge
University Press, 1957.

\bibitem{Da}
J. S. Dani,
{\it Properties of orbits under discrete groups },
J. Indian Math. Soc. (N.S.) {\bf 39} (1975),  189--217.

\bibitem{DaNo}
S. G. Dani  and A. Nogueira,
{\it On $SL(n,\bZ)_+$-orbits on $\bR^n$ and positive integral solutions of linear inequalities},
J. Number Theory (2009).

\bibitem{DaHe}
H. Davenport and H. Heilbronn, 
{\it Asymmetric inequalities for non-homogeneous linear forms}, 
J. London Math Soc. {\bf 22} (1947), 52--61.

\bibitem{Gu} 
A. Guilloux,
{\it A brief remark on orbits of {\rm SL(2,Z)} in the Euclidean plane},
Ergodic Theory and Dynamical Systems, published online 21 July 2009.

\bibitem{MaWe}
F. Maucourant and B. Weiss,
{\it Lattice actions on the plane revisited}, to appear, { \tt arXiv: 1001-4924v1}. 

\bibitem{No02}
 A. Nogueira,
{\it  Orbit distribution on $\bR^2$ under the natural action of $SL(2,\bZ)$},
Indag. Math.   {\bf 13} (2002), 103--124.

\bibitem{No03}
 A. Nogueira, 
{\it  Lattice orbit distribution on $\bR^2$ }, Ergodic Theory and Dynamical Systems, published online 21 July 2009.

\vfill\eject

\noindent Michel Laurent   \hfill{Arnaldo Nogueira}

\noindent Institut de Math\'ematiques de Luminy
\hfill{Institut de Math\'ematiques de Luminy}

\noindent C.N.R.S. -  U.M.R. 6206 - case 907
\hfill{C.N.R.S. -  U.M.R. 6206 - case 907}

\noindent  163, avenue de Luminy     \hfill{163, avenue de Luminy}

\noindent 13288 MARSEILLE CEDEX 9  (FRANCE)
\hfill{13288 MARSEILLE CEDEX 9  (FRANCE)}

\vskip2mm

\noindent {\tt nogueira@iml.univ-mrs.fr}
\hfill{{\tt laurent@iml.univ-mrs.fr}}

\endthebibliography

\bye